\documentclass{amsart}
\usepackage{amsmath, amstext, amsbsy, amssymb, amscd, graphics}

\newtheorem{theorem}{Theorem}[section]
\newtheorem{lemma}[theorem]{Lemma}
\newtheorem{proposition}[theorem]{Proposition}

\theoremstyle{definition}

\theoremstyle{remark}
\newtheorem{remark}[theorem]{Remark}
\newtheorem{conjecture}[theorem]{Conjecture}
\newtheorem{question}[theorem]{Question}

\numberwithin{equation}{section}

\newcommand{\B}{\mathbb B}

\newcommand{\C}{ \mathbb C }

\newcommand{\g}{{\gamma}}
\newcommand{\G}{{\Gamma}}
\newcommand{\Hn}{\mathcal{H}_n}

\newcommand{\la}{\lambda}

\newcommand{\Z}{ \mathbb Z }
\newcommand{\ZZ}{ \mathcal Z }
\def\N{\mathbb N}
\def\al{\alpha}

\def\H{{\mathcal H}}

{\vskip-\lastskip\medskip
  \noindent
  {\em #1.}\enspace
  }%
{\qed\par\medskip
  }

\begin{document}
\title[]
  {The centers of Iwahori-Hecke algebras are filtered}

\author[Andrew Francis]{Andrew Francis}
\address{School of Computing and Mathematics, University of Western Sydney, NSW 1797, Australia}
\email{a.francis@uws.edu.au}

\author[Weiqiang Wang]{Weiqiang Wang}
\address{Department of Mathematics, University of Virginia,
Charlottesville, VA 22904, U.S.A.} \email{ww9c@virginia.edu}


\begin{abstract}
We show that the center of the Iwahori--Hecke algebra of the
symmetric group $S_n$ carries a natural filtered algebra structure,
and that the structure constants of the associated graded algebra are
independent of $n$.
A series of conjectures and open problems are also included.
\end{abstract}

\maketitle
\date{}

\section{Introduction}

\subsection{The main results}

The class elements introduced by Geck-Rouquier \cite{GR} form a
basis for the center $\ZZ(\Hn)$ of the Iwahori-Hecke algebra $\Hn$
of type $A$ over the ring $\Z [\xi]$, where the indeterminate
$\xi$ is related to the familiar one $q$ by $\xi =q-q^{-1}.$ In
this paper, we shall parameterize these class elements $\G_\la(n)$
by partitions $\la$ satisfying $|\la| +\ell(\la) \le n$ which are
the so-called modified cycle types, just as Macdonald
\cite[pp.131]{Mac} does for the usual class sums of the symmetric
group $S_n$. Write the multiplication in $\mathcal Z(\H_n)$ as %%
\begin{equation} \label{eq:mult}
\G_\la(n) \G_\mu(n) =\sum_{\nu} k_{\la \mu}^\nu(n) \ \G_\nu(n).
\end{equation}

The main result of this Note is the following theorem on these
structure constants $k_{\la \mu}^\nu (n)$.
\begin{theorem} \label{th:main}
\begin{enumerate}
\item For any $n$, $k_{\la \mu}^\nu (n)$ is a polynomial in $\xi$
with non-negative integral coefficients. Moreover, $k_{\la
\mu}^\nu (n)$ is an even (resp., odd) polynomial in $\xi$ if and
only if ${|\la| +|\mu| -|\nu|}$ is even (resp., odd).

\item We have $k_{\la \mu}^\nu (n) =0$ unless $|\nu| \leq |\la|
+|\mu|$.

\item If $|\nu| = |\la| +|\mu|$, then $k_{\la \mu}^\nu (n)$ is
independent of $n$.
\end{enumerate}
\end{theorem}
It follows from (2) and (3) that the center $\mathcal Z(\H_n)$ is
naturally a filtered algebra and the structure constants of the
associated graded algebra are independent of $n$. We in addition
formulate several conjectures, including
Conjecture~\ref{conj:str.cons.polys.in.n} which simply states that
$k_{\la \mu}^\nu (n)$ are polynomials in $n$, and further
implications on the algebra generators of $\mathcal Z(\H_n)$.

\subsection{Motivations and connections}

Our motivation comes from the original work of Farahat-Higman
\cite{FH} on the structures of the centers of the {\em integral}
symmetric group algebras $\Z S_n$. Indeed, Theorem~\ref{th:main}
specializes at $\xi = 0$ to some classical results of {\em loc.
cit.}. In addition, Conjecture~\ref{conj:str.cons.polys.in.n} in the
specialization at $\xi = 0$ (which is a theorem in {\em loc. cit.})
when combined with the specialization of Theorem~\ref{th:main}
allowed them to define a universal algebra $\mathcal K$ governing
the structures of the centers $\mathcal Z(\Z S_n)$ for all $n$
simultaneously. Farahat-Higman further developed this approach to
establish a distinguished set of generators for $\mathcal Z(\Z
S_n)$, which is now identified as the first $n$ elementary symmetric
polynomials in the Jucys-Murphy elements. This has applications to
blocks of modular representations of $S_n$.

This Note arose from the hope that the results of Farahat--Higman
might be generalized to the Iwahori-Hecke algebra setup and in particular
it would provide a new conceptual proof of the Dipper-James
conjecture. Recently, built on the earlier work of Mathas
\cite{Mat}, the first author and Graham \cite{FG} obtained a first
complete combinatorial proof of the Dipper-James conjecture that
the center of $\Hn$ is the set of symmetric polynomials in
Jucys--Murphy elements of $\Hn$ (cf. \cite{DJ} for an earlier
proof of a weaker version). A basic difficulty in completing the
Farahat-Higman approach for Iwahori-Hecke algebras is that no
compact explicit expression is known for the Geck-Rouquier
elements $\G_\la(n)$ (see however \cite{F} for a useful
characterization). Our Theorem~\ref{th:main} is a first positive
step along the new line.

In another direction, the results of Farahat-Higman 
have been partially generalized by the second author \cite{W} to the
centers of the group algebras of wreath products $G\wr S_n$ for an
arbitrary finite group $G$ (e.g. a cyclic group $\Z_r$), and  these
centers are closely related to the cohomology ring structures of
Hilbert schemes of points on the minimal resolutions. It will be
interesting to develop the Farahat-Higman type results for the
centers of the Iwahori-Hecke algebra of type $B$ or more generally
of the cyclotomic Hecke algebras which are $q$-deformation of the
group algebra $\Z (\Z_r \wr S_n)$.

\subsection{}
This Note is organized as follows. In Section~\ref{sec:proof}, we
prove the three parts of our main Theorem~\ref{th:main} in
Propositions~\ref{positivity}, \ref{str.cons.upper.tri} and
\ref{prop:str.cons.indep.n}, respectively. In
Section~\ref{sec:conj}, we formulate several conjectures, which
are Iwahori-Hecke algebra analogues of some results of Farahat-Higman. We
conclude this Note with a list of open questions.

\section{Proof of the main theorem} \label{sec:proof}

\subsection{The preliminaries}
\label{subsec:pre}

 Let $S_n$ be the symmetric group in $n$ letters
generated by the simple transpositions $s_i =(i,i+1)$,
$i=1,\dots,n-1$.

Let $\xi$ be an indeterminate. The {\em Iwahori-Hecke algebra
$\H_n$} is the unital $\Z [\xi]$-algebra generated by $T_i$ for
$i=1,\dots,n-1$, satisfying the relations
\begin{eqnarray}
T_iT_{i+1}T_i &=&T_{i+1}T_iT_{i+1}
 \nonumber \\
T_iT_j &=&T_jT_i, \quad |i-j|>1
 \nonumber \\
T_i^2 &=& 1+\xi T_i. \label{eq:order}
\end{eqnarray}
The  order relation (\ref{eq:order}) comes from the more familiar
$(T_i-q)(T_i+q^{-1})=0$ via the identification $\xi =q-q^{-1}.$
If $w=s_{i_1}\cdots s_{i_r} \in S_n$ is a reduced expression
(where $r$ will be referred to as the {\em length} of $w$ in this
case), then define $T_w:=T_{i_1}\cdots T_{i_r}$. It is well known
that the Iwahori-Hecke algebra is a free $\Z[\xi]$-module with
basis $\{T_w\mid w\in S_n \}$, and it is a deformation of the
integral group algebra $\Z S_n$.

The {\em Jucys-Murphy elements} $L_i$ $(1\le i \le n)$ of $\Hn$
are defined to be $L_1=0$ and
$$
L_i =\sum_{1\le k<i} T_{(k,i)},
$$
for $i\ge 2$.

Given $w \in S_n$ with cycle-type $\rho =(\rho_1,\ldots,
\rho_t,1,\ldots,1)$ for $\rho_i>1$, we define the \textit{modified cycle-type} of
$w$ to be $ \tilde{\rho}  =(\rho_1-1,\ldots, \rho_t-1)$, following
Macdonald \cite[pp.131]{Mac}. Given a partition $\la$, let
$C_{\la}(n)$ denote the conjugacy class of $S_n$ containing all
elements of modified type $\la$ if $|\la| + \ell(\la) \leq n$.
Accordingly,  let
$c_{\la}(n)$ denote the class sum of $C_\la(n)$ if $|\la| +
\ell(\la) \leq n$, and denote $c_{\la}(n) =0$ otherwise.

The center $\ZZ(\Hn)$ of the Iwahori-Hecke algebra is free over
$\Z[\xi]$ of rank equal to the number of partitions of $n$; that
is, it has a basis indexed by the conjugacy classes of $S_n$ (see
\cite{GR}). In this paper, we shall parameterize these
Geck-Rouquier class elements $\G_\la(n)$ by the modified cycle
types $\la$. The elements $\G_\la(n)$ for $|\la| + \ell(\la) \leq
n$ are characterized by the following two properties among the
central elements of $\Hn$ \cite{F}:
\begin{enumerate}
\item[(i)] The $\Gamma_\lambda(n)$ specializes at $\xi=0$ to the
class sum $c_\lambda(n)$;

\item[(ii)] The difference $\Gamma_\lambda(n)-\sum_{w\in
C_\lambda(n)}T_w$ contains no minimal length elements of any
conjugacy class.
\end{enumerate}
In addition, we set $\G_{\la}(n) =0$ if $|\la| + \ell(\la) > n$.

\subsection{The structure constants as positive integral polynomials}

By inspection of the defining relations, $\Hn$ as a $\Z$-algebra
is $\Z_2$-graded by declaring that $\xi$ and $T_i$ $(1\le i \le
n-1)$ have $\Z_2$-degree (or parity) $1$ and each integer has
$\Z_2$-degree $0$.

\begin{lemma}\label{lem:Z2}
Every $\G_\la(n)$  is homogeneous in the above $\Z_2$-grading with
$\Z_2$-degree equal to $|\la| \mod 2$.
\end{lemma}
\begin{proof}

There is a constructive algorithm \cite[pp.14]{F} for producing
the elements $\G_\lambda(n)$.
This finite algorithm begins with the sum of $T_w$ with minimal
length elements $w$ from the conjugacy class $C_\la(n)$, then at
each repeat of this algorithm, the only additions are of form (i)
$T_w \to T_w+T_{s_iws_i}$, (ii) $T_{s_iw}\to T_{ws_i}+T_{s_iw}+\xi
T_{s_iws_i}$, or (iii) $T_{ws_i}\to T_{ws_i}+T_{s_iw}+\xi
T_{s_iws_i}$, and the algorithm eventually ends up with the
element $\G_\la(n)$. Each of the three type of additions clearly
preserves the $\Z_2$-degree. As the minimal length elements have
the same parity as $|\lambda|$, this proves the lemma.
\end{proof}

Denote by $\N$ the set of non-negative integers.

\begin{proposition} \label{positivity}
For any given $n$, $k_{\la \mu}^\nu (n)$ is a polynomial in $\xi$
with non-negative integral coefficients. Moreover, $k_{\la
\mu}^\nu (n)$ is an even (respectively, odd) polynomial in $\xi$
if and only if ${|\la| +|\mu| -|\nu|}$ is even (respectively,
odd).
\end{proposition}

\begin{proof}
As is seen in the proof of Lemma~\ref{lem:Z2}, the class elements
are in the positive cone $\mathcal Z(\Hn)^+:=\mathcal
Z(\Hn)\cap\sum_{w\in S_n}\N[\xi]T_w$. Because of the positive
coefficients in the order relation (\ref{eq:order}), the positive
cone is closed under additions and products. Since the class
elements contain minimal length elements from exactly one conjugacy
class and contains those minimal length elements with coefficient
$1$ \cite{F} (see Sect.~\ref{subsec:pre} above), the coefficient of
a class element in a central element $C$ is precisely the
coefficient of the corresponding minimal length element in the $T_w$
expansion of $C$. This shows that $k_{\la \mu}^\nu (n)\in \N [\xi]$.

The more refined statement on when $k_{\la \mu}^\nu (n)$ is an
even or odd polynomial in $\xi$ follows now from
Lemma~\ref{lem:Z2}.
\end{proof}

\subsection{The filtered algebra structure on the center}

Let $m_\mu(n)$ be the monomial symmetric polynomial in the
(commutative) Jucys-Murphy elements $L_1, \ldots, L_n$,
parameterized by a partition $\mu$. It is known that $m_\mu(n) \in
\ZZ(\Hn)$. Some relations between the class elements and the
monomial symmetric polynomials in Jucys--Murphy elements are
summarized as follows (see \cite{Mat, FG, FJ}).

\begin{lemma}\label{lem:cl.elts.w.mons}
\begin{enumerate}
\item For any partition $\la$, we can express $m_\lambda (n)$ in
terms of the class elements $\Gamma_\mu (n)$ as
$$m_\lambda (n) =\sum_{|\mu| \le |\la|} b_{\la\mu}(n) \  \G_\mu (n)
\ \quad \text{for } b_{\la\mu}(n) \in \Z[\xi].
$$
(We set $b_{\la\mu}(n) =0$ if $|\mu| +\ell(\mu) >n$, or equivalently
if $\G_{\mu}(n) =0$.)
 \label{lem:mons.to.cl.elts}

\item For any $\la$, the coefficients $b_{\la\mu}(n)$ with $|\mu|=|\lambda|$ are
independent of $n$. \label{indep-n}

\item Let $\la$ be a partition with $|\la| +\ell(\la) \leq n$.
Then, each $\Gamma_\lambda (n)$ is equal to $m_\lambda (n)$ plus a
$\Z[\xi]$-linear combination of $m_\mu (n)$ with
$|\mu|<|\lambda|$.\label{lem:cl.elts.to.mons}
\end{enumerate}
\end{lemma}
\begin{proof}
Part \eqref{lem:mons.to.cl.elts} is a consequence
of~\cite[Theorems~2.7 and 2.26]{Mat}.

By \cite[Lemma 5.2]{FG}, the coefficient of $T_w$ (for an
increasing $w\in S_n$ of the right length) in a so-called
quasi-symmetric polynomial in Jucys--Murphy elements is
independent of $n$. Monomial symmetric polynomials are just sums
of the corresponding quasi-symmetric polynomials, independently of
$n$. This proves \eqref{indep-n}.

Part \eqref{lem:cl.elts.to.mons} can be read off from the proof
of~\cite[Theorem 4.1]{FJ}.
\end{proof}

\begin{proposition}\label{str.cons.upper.tri}
We have  $k_{\la \mu}^\nu (n) =0$ unless $|\nu| \leq |\la|
+|\mu|$.
\end{proposition}

\begin{proof}
By Lemma~\ref{lem:cl.elts.w.mons}~\eqref{lem:cl.elts.to.mons}, the
product $\Gamma_\lambda (n)\Gamma_\mu (n)$ is equal to $m_\lambda
(n) m_\mu (n)$ plus a linear combination of products
$m_{\lambda'}(n) m_{\mu'}(n)$ satisfying
$|\lambda'|+|\mu'|<|\lambda|+|\mu|$. A product of monomial
symmetric polynomials $m_\alpha (n) m_\beta (n)$ is a sum of
monomial symmetric polynomials $m_\gamma (n)$ satisfying
$|\alpha|+|\beta|=|\gamma|$. Consequently, $\Gamma_\lambda(n)
\Gamma_\mu(n)$ is a sum of $m_\gamma (n)$ with partitions $\g$ of
size at most $|\lambda|+|\mu|$. The proposition follows now by
Lemma~\ref{lem:cl.elts.w.mons}~\eqref{lem:mons.to.cl.elts}.
\end{proof}

Assign degree $|\la|$ to the basis element $\G_\la(n)$ and let
$\ZZ (\Hn)_{m}$ be the $\Z[\xi]$-span of $\G_\la(n)$ of degree at
most $m$, for $m\ge 0$. Then, Proposition~\ref{str.cons.upper.tri}
provides a filtered algebra structure on the center $\ZZ (\Hn)
=\cup_m \ZZ (\Hn)_{m}$.

\begin{remark}
Denote by $A_n$ the (commutative) $\Z[\xi]$-subalgebra of $\H_n$
generated by the Jucys-Murphy elements $L_1,\ldots, L_n$. The
algebra $A_n$ is filtered by the subspaces $A_n^{(m)}$ $(m\ge 0)$
spanned by all products $L_{i_1}\cdots L_{i_m}$ of $m$
Jucys-Murphy elements.
The filtrations on $A_n$ and $\ZZ(\H_n)$ are compatible with each
other by the inclusion $\ZZ(\Hn) \subset A_n$. However the algebra
$\Hn$ does not seem to admit a natural filtration which is
compatible with the one on $A_n$ by inclusion $A_n \subset \Hn$.
This is very different from the symmetric group algebra $\Z S_n$,
which admits such a filtration by assigning degree $1$ to every
transposition $(i,j)$.
\end{remark}

\subsection{The graded algebra $\text{gr} \ZZ(\Hn)$}

\begin{proposition}\label{prop:str.cons.indep.n}
If $|\nu| = |\la| +|\mu|$, then $k_{\la \mu}^\nu (n)$ is
independent of $n$.

(In this case, we shall write $k_{\la \mu}^\nu (n)$ as $k_{\la
\mu}^\nu.$)
\end{proposition}

\begin{proof}
By the definition of $k_{\la \mu}^\nu (n)$, we can assume without
loss of generality that $|\la| +\ell(\la) \leq n$ and $|\mu|
+\ell(\mu) \leq n$.

By Lemma~\ref{lem:cl.elts.w.mons}~\eqref{lem:cl.elts.to.mons},
$\Gamma_\lambda(n) \Gamma_\mu(n) =m_\lambda(n) m_\mu(n) +X$, where
$X$ is a linear combination of products of monomials whose
combined partition size is less than $|\lambda|+|\mu|$. The
monomials appearing in $X$ correspond to partitions of size less
than $|\nu|= |\lambda|+|\mu|$, and thus will not contribute to
$k_{\la \mu}^\nu (n)$ by
Lemma~\ref{lem:cl.elts.w.mons}~\eqref{lem:mons.to.cl.elts}. The
product $m_\lambda (n) m_\mu(n)$ is a sum of monomials $m_\al(n)$
satisfying $|\al|=|\lambda|+|\mu|$ with coefficients independent
of $n$; the contribution of each such $m_\al (n)$ to $k_{\la
\mu}^\nu (n)$ is independent of $n$ by
Lemma~\ref{lem:cl.elts.w.mons}~\eqref{indep-n}. Summing all these
contributions produces $k_{\la \mu}^\nu (n)$ which is independent
of $n$.
\end{proof}

Proposition~\ref{prop:str.cons.indep.n} is equivalent to the
statement that all the structure constants of the graded algebra
$\text{gr} \ZZ (\Hn)$ associated to the filtered algebra $\ZZ
(\Hn)$ are independent of $n$.

\subsection{Examples}
\label{example} We provide some explicit calculations of the
structure constants $k_{\la\mu}^\nu(n)$ for the multiplication
between $\G_\la(n)$, with $n=3, 4, 5$. For the sake of notational
simplicity, we will write $\G_\la(n)$ as $\G_\la$ with $n$ dropped
in the following examples.  We also drop parentheses in the subscripts of class elements, denoting $\Gamma_{(\lambda_1,\dots,\lambda_k)}$ by $\Gamma_{\lambda_1,\dots,\lambda_k}$.
The square brackets denote the
top-degree parts of each product. The compatibility of these
examples with Theorem~\ref{th:main} is manifest.

\vspace{.2cm}

\noindent (1) Let $n=3$. In $\ZZ (\mathcal H_3)$, we have
\[\Gamma_1\Gamma_1=\left[(\xi^2+3)\Gamma_2\right]+2\xi\Gamma_1+3\Gamma_\emptyset.
\]
\noindent (2) Let $n=4$. In $\ZZ (\mathcal H_4)$, we have
\begin{align*}
\Gamma_1\Gamma_1
&=\left[(\xi^2+3)\Gamma_2+(\xi^2+2)\Gamma_{1,1}\right]
 +3\xi\Gamma_1+6\Gamma_\emptyset, \\
\Gamma_1\Gamma_2 &=\left[(\xi^4+4\xi^2+4)\Gamma_3\right]
  +(2\xi^3+6\xi)\Gamma_2+(2\xi^3+4\xi)\Gamma_{1,1}\\
&\qquad +(3\xi^2+4)\Gamma_1+4\xi\Gamma_\emptyset, \\
\Gamma_1\Gamma_{1,1} &=\left[(\xi^2+2)\Gamma_3\right]
 +2\xi\Gamma_2+\xi\Gamma_{1,1} +\Gamma_1.
\end{align*}

\noindent (3) Let $n=5$. In $\ZZ (\mathcal H_5)$, we have
{\allowdisplaybreaks
\begin{align*}
\Gamma_1\Gamma_1
  &=\left[(\xi^2+3)\Gamma_2+(\xi^2+2)\Gamma_{1,1}\right]+4\xi\Gamma_1+10\Gamma_\emptyset\\
\Gamma_1\Gamma_2
  &=\left[(\xi^4+4\xi^2+4)\Gamma_3+(\xi^4+2\xi^2+1)\Gamma_{2,1}\right]\\
&\qquad +(3\xi^3+8\xi)\Gamma_2+(3\xi^3+4\xi)\Gamma_{1,1}
  +(6\xi^2+6)\Gamma_1+10\xi\Gamma_\emptyset\\
\Gamma_1\Gamma_{1,1}&=\left[(\xi^2+2)\Gamma_3+(2\xi^2+3)\Gamma_{2,1}\right]
  +2\xi\Gamma_2+4\xi\Gamma_{1,1} +3\Gamma_1\\
\Gamma_1\Gamma_3&=\left[(\xi^6+6\xi^4+10\xi^2+5)\Gamma_4\right]\\
&\qquad +(2\xi^5+10\xi^3+13\xi)\Gamma_3+(2\xi^5+8\xi^3+7\xi)\Gamma_{2,1}\\
&\qquad +(3\xi^4+10\xi^2+6)\Gamma_2+( 3\xi^4+8\xi^2+4)\Gamma_{1,1}
  +(4\xi^3+6\xi)\Gamma_1+5\xi^2\Gamma_\emptyset\\
\Gamma_2\Gamma_2&=\left[(\xi^8+7\xi^6+16\xi^4+15\xi^2+5)\Gamma_4\right]\\
&\qquad +(2\xi^7+14\xi^5+29\xi^3+19\xi)\Gamma_3+(2\xi^7+13\xi^5+22\xi^3+11\xi)\Gamma_{2,1}\\
&\qquad +(3\xi^6+20\xi^4+32\xi^2+7)\Gamma_2+(3\xi^6+19\xi^4+26\xi^2+8)\Gamma_{1,1}\\
&\qquad +(4\xi^5+25\xi^3+27\xi)\Gamma_1+(5\xi^4+30\xi^2+20)\Gamma_\emptyset\\
\Gamma_2\Gamma_{1,1}&=\left[(\xi^6+6\xi^4+10\xi^2+5)\Gamma_4\right]\\
&\qquad +(2\xi^5+10\xi^3+11\xi)\Gamma_3+(2\xi^5+9\xi^3+9\xi)\Gamma_{2,1}\\
&\qquad +(3\xi^4+11\xi^2+6)\Gamma_2+(3\xi^4+9\xi^2+4)\Gamma_{1,1}
 +(4\xi^3+7\xi)\Gamma_1+5\xi^2\Gamma_\emptyset.
\end{align*}
 }

\subsection{A universal graded algebra}

Introduce a graded $\Z[\xi]$-algebra $\mathcal G$ with a basis
given by the symbols $\G_\la$, where $\la$ runs over all
partitions, and with multiplication given by
$$
\G_\la \G_\mu =\sum_{|\nu| =|\la|+|\mu|} k_{\la \mu}^\nu \G_\nu.
$$
By Propositions~\ref{positivity} and \ref{prop:str.cons.indep.n},
the structure constants $k_{\la \mu}^\nu$ are independent of $n$
and actually lie in $\N [\xi^2]$. Furthermore, we have surjective
homomorphisms  $\mathcal G \rightarrow \text{gr} \ZZ [\Hn]$ for
all $n$, which send each $\G_\la$ to $\G_\la(n)$. The following
proposition is immediate.

\begin{proposition}
The $\Z[\xi]$-algebra $\mathcal G$ is commutative and associative.
\end{proposition}

Below for the one-row partition $(m)$, we shall write $\G_{(m)}$
simply as $\G_m$.

\begin{theorem}
The $\mathbb Q(\xi)$-algebra $\mathbb
Q(\xi)\otimes_{\Z[\xi]}\mathcal G$ is a polynomial algebra with
generators $\G_{m}$, $m=1,2,\ldots$.
\end{theorem}

\begin{proof}
Given a partition $\la =(\la_1,\la_2,\ldots)$,  the product is of
the form
$$\G_{\la_1}\G_{\la_2}\cdots =\sum_\mu d_{\la\mu}(\xi)\  \G_\mu
$$
for $d_{\la\mu}(\xi) \in \N[\xi].$ As $\xi$ goes to $0$, $\G_\mu$
goes to the class sum $c_\mu$ and $d_{\la\mu}(\xi)$ specifies to
the structure constant $d_{\la\mu}$ as defined in
\cite[pp.132]{Mac} which we recall:
$$c_{\la_1}c_{\la_2}\cdots =\sum_\mu d_{\la\mu}c_\mu.$$
It is known therein that the (integral) matrix $[d_{\la\mu}]$ for
$|\la| =|\mu| =k$ with any $k$ is triangular with respect to the
dominance ordering of partitions and all its diagonal entries are
nonzero, thus the matrix $[d_{\la\mu}]$ is invertible over
$\mathbb Q$. This forces the matrix $[d_{\la\mu} (\xi)]$
invertible over the field $\mathbb Q(\xi)$. Thus each $\G_\mu$ is
generated by $\G_1, \G_2, \ldots$ over $\mathbb Q(\xi)$. By
definition, the elements $\G_\mu$ for all partitions $\mu$ are
linearly independent. Thus the theorem follows by comparing the
graded dimensions of the algebra $\mathbb Q(\xi)
\otimes_{\Z[\xi]}\mathcal G$ and the polynomial algebra in $\G_m$,
$m=1,2,\cdots$.
\end{proof}

It is not clear though if the matrix  $[d_{\la\mu} (\xi)]$ remains
triangular. A similar type of phenomenon with a negative answer
appears in the example of $\widetilde{M}_3$ in Mathas
\cite[p.310]{Mat}.

\section{Conjectures and discussions}
\label{sec:conj}

\subsection{Several conjectures}

We expect the following conjecture to hold.
\begin{conjecture}\label{conj:str.cons.polys.in.n}
Given partitions $\la, \mu$ and $\nu$, there exists a polynomial
$f_{\la \mu}^\nu$ in one variable with coefficients in $\mathbb
Q[\xi]$, such that $f_{\la \mu}^\nu (n) =k_{\la \mu}^\nu (n)$ for
all $n$.
\end{conjecture}

Recall $b_{\la\mu}(n)$ from
Lemma~\ref{lem:cl.elts.w.mons}~\eqref{lem:mons.to.cl.elts}.
Similarly, we conjecture that there exists polynomials  $g_{\la \mu}
(x)$ with coefficients in $\mathbb Q[\xi]$, such that $g_{\la
\mu}(n) =b_{\la \mu} (n)$ for all $n$.

The specialization  at $\xi=0$ of
Conjecture~\ref{conj:str.cons.polys.in.n} is a result in
\cite{FH}. {\em Below we shall assume that
Conjecture~\ref{conj:str.cons.polys.in.n} holds.}

Set $\B$ to be the ring of polynomials in $\mathbb Q[x]$ which
take integer values at integers.  We can define a
$\B[\xi]$-algebra $\mathcal K$ with a basis given by the symbols
$\G_\la$, where $\la$ runs over all partitions, and the
multiplication given by
$$
 \G_\la \G_\mu =\sum_\nu f_{\la \mu}^\nu  \G_\nu.
$$
Since $f_{\la \mu}^\nu =0$ unless $|\nu| \leq |\la| +|\mu|$, the
algebra $\mathcal K$ is an algebra filtered by $\mathcal K_m$
($m\ge 0$) which is the $\B[\xi]$-span of $\G_{\la}$ with $|\la|
\le m$. We have a natural surjective homomorphism of filtered
algebras
$$
p_n: \mathcal K   \longrightarrow \ZZ(\H_n)
$$
given by
\begin{align*}
p_n \big (\sum f_\la \G_\la \big ) & = \sum_\la f_\la(n) \
\G_\la(n).
\end{align*}
The algebra $\mathcal G$ introduced earlier becomes the associated
graded algebra for the filtered algebra $\mathcal K$ up to a base
ring change, i.e., $\text{gr} \mathcal{K}
=\B[\xi]\otimes_{\Z[\xi]} \mathcal G$. For $r \ge 1$, set
$$E_r :=\sum_{|\la|=r} \G_\la.$$

\begin{conjecture}\label{conj:algK}
The $\B[\xi]$-algebra $\mathcal K$ is generated by $E_r$,
$r=1,2,\cdots$.
\end{conjecture}
Conjecture~\ref{conj:algK} in the specialization with $\xi=0$ is a
main theorem of Farahat-Higman \cite{FH}.

Note (cf. e.g. \cite{FG}) that
$$
E_r(n) :=\sum_{|\la|=r} \G_\la(n) \in \ZZ (\Hn)
$$
can be interpreted as the $r$-th elementary symmetric function
in the $n$ Jucys--Murphy elements. If Conjecture~\ref{conj:algK}
holds, then the surjectivity of the homomorphism $p_n: \mathcal
K\rightarrow \mathcal Z(\H_n)$ implies that the center $\mathcal
Z(\H_n)$ is generated by $E_r(n)$, $1\le r \le n-1$.  That would
provide a new and conceptual proof of the Dipper-James
conjecture~\cite{DJ, FG} along with additional results of
independent interest.

\subsection{Open questions and discussions}

A fundamental difficulty in pursuing the approach of \cite{FH} for
the Iwahori-Hecke algebras is present in the following.
\begin{question}
Find an explicit expression for the elements $\G_\la(n)$ for all
$\la$.
\end{question}

The more challenging Conjecture~\ref{conj:algK} is likely to follow
from an affirmative answer to Question~\ref{Qcal} below, as a
similar calculation in the case of symmetric groups plays a key role
in the original approach of \cite{FH}.
\begin{question} \label{Qcal}
Calculate the structure constants $k_{\la\, (m)}^\nu$ with $|\nu|
=|\la| +m$.
\end{question}

Recall the positivity and integrality from
Theorem~\ref{th:main}~(1).

\begin{question}
Find a combinatorial or geometric interpretation of the positivity
and integrality of the structure constants $k_{\la\mu}^\nu(n)$ as
polynomials in $\xi$.
\end{question}

With the connections between results of Farahat-Higman and
cohomology rings of Hilbert schemes of $n$ points on the affine
plane in mind (cf. \cite{W}), we post the following.

\begin{question} \label{QHilbertC2}
Are there any connections between the center $\ZZ (\Hn)$ and the
equivariant $K$-group of Hilbert schemes of $n$ points on the affine
plane?
\end{question}

Let $W$ be an arbitrary finite Coxeter group. The group algebra
$\Z W$ is naturally filtered by assigning degree $1$ to each
reflection (not just simple reflection) and degree $r$ to any
element $w\in W$ with a reduced expression in terms of reflections
of minimal length $r$. This induces a filtered algebra structure
on the center $\ZZ (\Z W)$ as elements of a conjugacy class have
the same degree. In the cases of types $A$ and $B$, this
definition agrees with the notion of degree for general wreath
products introduced in~\cite{W}. (In the case of symmetric groups,
the degree of $c_\la(n)$ coincides with $|\la|$.) The
Geck-Rouquier basis has been defined for centers of the integral
Iwahori-Hecke algebras $\mathcal H_W$ associated to any such $W$
\cite{GR}, and its characterization (as in Sect.~\ref{subsec:pre})
holds in this generality \cite{F}. Note that the generalization of
Theorem~\ref{th:main}~(1) to all such $\mathcal H_W$ holds with
the same proof. We ask for a generalization of
Theorem~\ref{th:main}~(2) as follows.

\begin{question}  \label{QtypeB}
Let $W$ be an arbitrary finite Coxeter group. If we apply the
notion of degree above to the Geck-Rouquier elements in the center
$\ZZ (\mathcal H_W)$, does it provide an algebra filtration on
$\ZZ (\mathcal H_W)$?
\end{question}

We expect that the answer to Question~\ref{QtypeB}, at least for
Iwahori-Hecke algebras of type $B$, is positive. More generally, we
ask the following.
\begin{question}
Establish and characterize an appropriate basis of class elements
for the centers of the integral cyclotomic Hecke algebras associated
to the complex reflection groups $\Z_r \wr S_n$. Furthermore,
generalize the results of this Note and \cite{W} to the cyclotomic
Hecke algebra setup. Are there any connections between these centers
and equivariant $K$-groups of Hilbert schemes of points on the
minimal resolution $\widetilde{\C^2/\Z_r}$?
\end{question}
It will be already nontrivial and of considerable interest to answer
the question for the Iwahori-Hecke algebras of type $B$
corresponding to $r=2$.

\vspace{.3cm}

{\bf Acknowledgment.} The authors thank the organizers for the
high level International Conference on Representation Theory in
Lhasa, Tibet, China in July 2007, where this collaboration was
initiated. W.W. is partially supported by an NSF grant.

\end{document}